\documentclass[11pt]{article}
 \usepackage{amsfonts}
 \usepackage{amssymb}

\newtheorem{deff}{Definition}

\newtheorem{propos}{Proposition}

\newtheorem{theorem}{Theorem}
\newtheorem{lemma}{Lemma}

\newcommand{\proof}{{\bf Proof.}~}

\newcommand\bqa {\begin{eqnarray}}
\newcommand\eqa {\end{eqnarray}}
\newcommand{\beq}{\begin{eqnarray}}
\newcommand{\eeq}{\end{eqnarray}}
\newcommand{\be}{\begin{array}}
\newcommand{\ee}{\end{array}}

 \newcommand{\pr}{\partial}

 \newcommand{\E}{{\cal E}}
 \newcommand{\s}{{\cal S}}
 \newcommand{\W}{{\cal W}}
 \newcommand{\V}{{\cal V}}

\newcommand{\C}{{\mathbb C}}
\newcommand{\R}{{\mathbb R}}
\newcommand{\Z}{{\mathbb Z}}
\newcommand{\Q}{{\mathbb Q}}

\newcommand{\HH}{{\mathbb H}}

\newcommand{\RR}{{\mathbb R}}

\newcommand{\A}{{\mathbb A}}

 \newcommand{\g}{{\mathfrak g}}
 \newcommand{\su}{{\mathfrak su}}
 \newcommand{\so}{{\mathfrak so}}
 \newcommand{\sll}{{\mathfrak sl}}

\begin{document}

   \def\sp{\mathfrak sp}
   \def\S{{\cal S}}

\begin{flushright}
 ITEP-TH 83-2000 \\
\end{flushright}

\vskip 10mm

\begin{center}
{ \bf \large Transgression on Hyperk\"{a}hler Manifolds and \\
   Generalized Higher Torsion Forms.}
\end{center}

\smallskip

\begin{center}
A.Gerasimov, A.Kotov
\end{center}

\centerline{\sl ITEP, B. Cheriomushkinskaya 25, Moscow, 117259 Russia}
\vskip 10mm





Transgression of the characteristic classes taking values in the
differential forms is a reach source of the interesting algebraic
objects. The examples include Chern-Simons and Bott-Chern forms
which are given by the transgression of the Chern character form.
Chern-Simons forms are defined for a vector bundle over an
arbitrary real manifold and are connected with the representation
of  combinations of Chern classes by the exact form
$\omega=d\phi$.  Bott-Chern forms  are defined for  holomorphic
hermitian vector bundles over
 K\"ahler  manifolds. These additional structures allow to use the
double transgression $\omega=\partial \overline{\partial}\phi $ to
define this invariant. Basically the existence of this
representation is a consequence of the action of the
multiplicative group of complex numbers $\mathbb{C}^*$ on the
cohomology of an arbitrary  K\"{a}hler manifold.

It is natural to guess that in the case when there is a bigger
group acting on the cohomology one should look for more involved
objects associated with vector bundles. In this paper we consider
the case of the action of the multiplicative group of quaternions
$\mathbb{H}^*$ on the cotangent bundle which induces the action of
$\mathbb{H}^*$ on the cohomology of the manifold. Supplying the
manifold with a metric compatible with the action of
$\mathbb{H}^*$ we get a hyperk\"{a}hler manifold. We propose a new
invariant of a hyperholomorphic bundle over a hyperk\"{a}hler
manifold connected with the Chern character form by the {\it
fourth} order "transgression"
$\omega=dd_Id_Jd_K \phi$. It takes
values in differential forms and its zero degree part is
hyprholomorphic analog of the logarithm of the holomorphic torsion
(holomorphic torsion is trivial for hyperk\"ahler manifolds). This
new hypertorsion seems first have appeared in the physical
literature \cite{BV}.

The  expression for the hypertorsion  in terms of the integration
over quaternionic projective plane proposed in this paper is a
direct generalization of the formula for the double transgression
\cite{GS}. The double transgression of the Chern character form in
terms of the integration over complex projective plane provides
the first example of the series of the regulator maps in algebraic
K-theory.
 We believe that the results of this paper imply (among other interesting
  applications) that
there is  a generalization of the regulator maps in algebraic
K-theory with the basic simplex being the configuration of linear
subspaces in the quaternionic linear spaces.

The paper is organized as follows. In the first part  we propose
the generalization of the Hodge $dd_c$-lemma for  compact
hyperk\"{a}hler manifolds. This leads to the fourth order
transgression of the differential forms. In the second part we
consider the transgression of the Chern classes of
hyperholomorphic bundles. Application of  the results of the first
section  gives the global construction of the fourth order
transgression of the Chern character of hyperholomorphic bundles.
Then we give the explicit local construction of this new invariant
for the important example of the infinite dimensional bundle
arising in the discussion  of the local families index theorem. We
define the higher analytic hypertorsion for families of
hyperholomorphic bundles on   compact hyperk\"ahler manifolds. An
explicit formula for the zero-degree part is given in terms of the
Laplace operators acting on sections of the vector bundle twisted
by the bundle of the differential forms.

{\bf Acknowledgements}: During the course of this work, the
authors benefited from helpful conversations with A.Levin. We also
grateful V.Rubtsov for useful comments. The work of A.G. was
partially supported by RFBR grant 98-01-00328 and Grant for the
Support of Scientific Schools 00-15-96557. The work of A.K. was
partially supported by RFBR grant 98-01-00327 and Grant for the
Support of Scientific Schools 00-15-99296.

\begin{center}
\large \bf Generalization of $dd_c$-lemma for hyperk\"ahler
manifolds
\end{center}

\vskip 3mm

To put the result of this section in the right perspective we
begin with the well-known cases of the Riemannian and K\"ahler
manifolds and then consider the case of the Hyperk\"ahler
manifold.

 Let us given a one-dimensional family of the closed
differential forms $\omega$ on the compact Riemannian manifold $M$
with the constant image in the  de Rham cohomology group.
 \bqa \label{real1}
 \omega(t) \in \Omega^{closed} \eqa
 \bqa \label{real2}
 \left[\delta_t \omega (t)\right]=0 \mbox{\, in \,} H^{\bullet}(M)
\eqa

  It implies that the variation of the differential form is exact:
  \bqa \delta_t \omega(t)=d\phi \eqa

  In the presence of the metric on the manifold $M$ one could chose
  the unique representative for the form $\phi$. Let $d$ be the  exterior
  derivative operator, $d^*$ be a conjugated operator with respect to the
  natural scalar product on the differential forms defined by the
  metric. Then  the corresponding Laplace operator $\Delta=dd^*+d^*d$
  and its  Green function $G$ are defined:
\bqa
 1=\mathcal{H}+\Delta G \eqa
 where $\mathcal{H}$ is the projector on the harmonic forms.
   Using the standard
  considerations \cite{GH}
  one gets the following explicit expression for the form $\phi$.
  \bqa \phi =d^*G\delta_t \omega(t) \eqa

 Consider now the case of the K\"ahler manifold $M$.
Let $I$ be an automorphism of the cotangent bundle corresponding
to integrable covariantly constant complex structure on $M$. One
could extended the action of $I$  on $k-$forms for arbitrary $k$
as: \bqa
 ad_I (\alpha_1 \wedge ...\wedge\alpha_k)=\sum_i
\alpha_1...\wedge I(\alpha_i)\wedge...\alpha_k. \eqa

\begin{deff}
{The Lie algebra $\g_M$ generated by $ad_I$ is called an isotropy
algebra, and the corresponding Lie group $G_I=U(1)$ is called an
isotropy group.}
\end{deff}
 Considering  $I$ as an element of the isotropy group $G_M$
 we have
 \bqa
 I (\alpha_1 \wedge
...\wedge\alpha_k)=I(\alpha_1)...\wedge
I(\alpha_i)\wedge...I(\alpha_k)\eqa

Given a complex structure $I$, one can introduce the differential
operator $d_I=[ad_I,d]$. Since $ad_I(\omega)=i(p-q)\omega$ for
$\omega$ of type $(p,q)$ we have $d_I=i(\partial_I
-\bar{\partial}_I)$. Considering $I$ as an element of the group
 ($I(\omega)=i^{p-q}\omega$), the new differential may be
represented as $d_I=IdI^{-1}$. These two equivalent
representations of $d_I$ immediately imply \bqa \label{basic}
d_I^2=dd_I+d_Id=0. \eqa

It is useful for further generalizations to introduce the
differential operator $d_x=x^0d+x^1d_I$ parameterized by the point
of the complex plane $(x^0+x^1i) \in \mathbb{C}$. We also define
the operator $\widehat{x}=x^0N+x^1ad_I $ where $N$ is the grading
operator acting as $k$ on the differential $k$-form. These
operators have the following obvious properties: \bqa
(i)&[ad_{\widehat{x}}, d_{y}]&=d_{xy} \label{relone} \\
(ii)&\{d_x,d_{y}\}&=d_x d_{y} +d_{y}d_x=0 \label{reltwo}\\
(iii)&\{d_x, d_{y}^* \}&=d_x d_{y}^* +d_{y}^*
d_x=Re(\bar{x}y)\Delta \label{relthree} \eqa

Let us given a one-dimensional family of the closed $G_M=U(1)$
invariant differential forms $\omega(t)$ on a compact K\"ahler
manifold $M$. Suppose this family has the constant image in de
Rham cohomology of $M$. Thus we have the conditions:

 \bqa \label{complex1} \omega(t) \in \Omega^{closed} \eqa
 \bqa \label{complex2} [\delta_t \omega(t)]=0 \mbox{\, in \,} H^{\bullet}(M)
 \eqa
 \bqa \label{complex3} [ad_I,\omega]=0 \eqa

In particular the last condition implies that the form $\omega$ is
$d_I$-closed and thus  the following theorem (Hodge
"$dd_c$-lemma") is applied:
\begin{theorem}\cite{GH}
Let $\omega$ be a $d$-exact and $d_I$-closed form on a compact
K\"ahler manifold. Then: $\omega=dd_I \chi$
\end{theorem}
Therefore the variation of the form $\omega$ is $dd_I$-exact:
 \bqa \label{comprep} \delta_t \omega(t)=dd_I \chi \eqa
 The explicit formula for $\chi$ may be given in terms of the
 Laplace operator and its Green function (see e.g. \cite{GH}):

$$ \chi=dd^*Gd_Id_I^*G\delta \omega(t)= dd_Id^*d^*_IG^2 \delta_t
\omega(t) $$

Here we have used the relations (\ref{reltwo})(\ref{relthree}).
Note that the condition (\ref{complex3}) may be substituted by a
more strong condition on the complex valued differential form
$\omega$ to be an eigenvalue of the operator $ad_I$. This gives
rise to the same representation (\ref{comprep}).

The next case to consider is the differential forms on
hyperk\"ahler manifolds.

\begin{deff}\cite{Hit}
{A hyperk\"ahler manifold is a Riemannian manifold $M$ with three
complex structures $I$, $J$ and $K$ which satisfy the quaternionic
identities $I^2=J^2=K^2=IJK=-1$, such that $M$ is K\"ahler with
respect to any of these structures. }
\end{deff}

On a hyperk\"ahler manifold we have a family of integrable complex
structures $$C=\sigma^1 I+\sigma^2 J+\sigma^3 K$$
$$C^2=-(\sigma^1)^2-(\sigma^2)^2-(\sigma^3)^2=-1 $$
 parameterized by the points of the sphere $S^2$, such that
 $M$ is K\"ahler with respect to
any $C$.

Obviously $\g_M=\su (2)$, $G_M=SU(2)=Sp(1)$. This means that $$
[ad_I, ad_J]=2ad_K, \hspace{3mm} [ad_J, ad_K]=2ad_I, \hspace{3mm}
[ad_K, ad_I]=2ad_J. $$

Consider the differential operator $d_x=x^0 d +x^1 d_I +x^2 d_J
+x^3 d_K $ parameterized by the points on the quaternioinic plane
$x=x^0 +x^1 i+x^2 j+x^3 k\in \HH$ and the operator $\widehat{x}
=x^0N +x^1 ad_I+x^2 ad_J+x^3 adK$ acting on the differential forms.
These operators satisfy the following relations:
\begin{propos} \label{trgr1} {
\bqa (i)& [ad_{\widehat{x}}, d_{y}]&=d_{xy} \label{bf1}\\ (ii)&
\{d_x, d_{y}\}&=d_x d_{y} +d_{y}d_x=0 \label{bf2}\\ (iii)& \{d_x,
d_{y}^* \}&=d_x d_{y}^* +d_{y}^*
d_x=Re(\bar{x}y)\Delta,\label{bf3} \eqa }
\end{propos}
\proof (i) First let us prove that $[ad_{\widehat{x}},d]=d_{x}$.
Since the Levi-Chivita connection on $M$ is torsion-free, one has
the relation $d(\alpha)=d\xi^i \wedge\nabla_i(\alpha)$, where
$\xi^i$ are local coordinates. Operators $N$, $I$, $J$, $K$ are
covariantly constant and  we get $$ [ad_{\widehat{x}},d]=(x^0 +
x^1 I+x^2 J +x^3 K)(d\xi^i)\wedge \nabla_i=$$ $$= x^0 d+x^1
[ad_I,d] +x^2 [ad_J,d] +x^3 [ad_K,d]=d_x $$ Now we have  $$
[ad_{\widehat{x}} ,d_{y}]=(x^1I+x^2J+x^3K )(y^0 + y^1 I+y^2 J +y^3
K) (d\xi^i)\wedge \nabla_i = [ad_{\widehat{xy}},d]=d_{xy} $$ Using
these identities, one can infer that the action of the  invertable
quaternion $U$ by the conjugation on $d_x$ has a simple form: \bqa
U d_x U^{-1}=d_{U x}\eqa (ii) Obviously  $\{d_x, d_{y}\}$ is
conjugated to $\{d, d_{\bar{x}y}\}$ and  therefore (ii) is the
consequence of (\ref{basic}).(iii) Since the action of the
isotropy group $G_M$ is unitary one has $$Ud_x^*
U^{-1}=d_{Ux}^*.$$ Let $L^C=\omega_C\wedge $ be the Hodge operator
of the multiplication on the K\"ahler form for the complex
structure $C$ and $\Lambda_C=L_C^*$ be its conjugate. Kodaira's
identities $$ d_C^*=[\Lambda_C, d], \hspace{3mm}
d^*=-[\Lambda_C,d_C], $$ imply the following relations: $$ \{ d,
d_C^*\}=0,\hspace{3mm} \{ d^*, d_C\}=0  $$ $$
\Delta_{d}=\Delta_{d_C}, \hspace{3mm} \Delta_{d_C}=d_C d_C^*
+d_C^* d_C $$ From these identities and the formula for the
conjugation of $d_x$ and $d_x^*$, we easily derive:  $$\{d_x,
d_{y}^* \}=U\{d, d_{\bar{x}y}^*\}U^{-1}=Re(\bar{x}y)U\Delta
U^{-1}, $$ where $U=|x|^{-1}x$. Finally  (iii) follows from the
fact that $\Delta$ is $G_M-$invariant. $\blacksquare$

Consider a one-dimensional family of the closed
$G_M=Sp(1)$-invariant differential forms $\omega (t)$ on a compact
hyperk\"ahler manifold $M$ and suppose $\omega(t)$ has the
constant image in the de Rham cohomology:

 \bqa \label{quater1} \omega(t) \in \Omega^{closed} \eqa
 \bqa \label{quater2} [\delta_t \omega(t)]=0 \mbox{\, in \,} H^{\bullet}(M)
 \eqa
 \bqa \label{quater3} [ad_{G_M},\omega]=0 \eqa

The last condition implies that $\omega(t)$ is $d_C$-closed for
any C. We would like to show that the variation of the
differential form admits the following representation: \bqa
 \delta_t\omega=d d_I d_J d_K (\tau) \eqa
The variation of the form is exact and the result follows from the
proposition:
\begin{propos}
Let $\omega$ be an $d$-exact and $d_C$-closed  differential form
of order $k$ for any compatible complex structure $C$ on the
compact hyperk\"ahler manifold $M$. Then there exists a form
$\tau$ of order $k-4$ such that $\omega=d d_I d_J d_K (\tau)$.
\end{propos}
\proof Note that if $\omega$ is exact and $d_C$-closed then by
Hodge theorem  $\omega$ is $d_C-$exact for all $C$. In particular
$\omega$ is of type $(p,p)$ with respect to any $C$ from the
hyperk\"{a}hler family.

Let $G$ be the Green operator associated with the Laplacian
$\Delta$. Then by taking into account that $\omega$ is $d_C-$exact
we obtain $$ \omega=d_Cd_C^* G\omega. $$ From the relations
(\ref{bf1})-(\ref{bf3}) we have $$ \omega=(\prod_{C=I,J,K}
d_Cd_C^* G)\omega= d d_I d_J d_K (d^* d_I^* d_J^* d_K^* G^4
\omega). $$ Therefore taking $\tau=d^* d_I^* d_J^* d_K^* G^4
\omega$ we immediately obtain the desired formula $\omega=d d_I
d_J d_K (\tau)$. $\blacksquare$

For a four dimensional manifold this relation may be further
simplified.

\begin{propos} \label{lapl}
{Let $M$ be a hyperk\"{a}hler manifold of dimension $4$ and let
$\varphi$ be a smooth function with a compact support.
 Then $dd_Id_Jd_K(\varphi)=16 vol_M \Delta^2 (\varphi)$.}
\end{propos}
\proof

Taking into account the properties of the differentials and
K\"ahler forms under conjugation:
 \beq \label{ko2} \left\{
\begin{array}{ccc}
d_J &=& JdJ^{-1} \\ d_K &=& -Jd_IJ^{-1}
\end{array}
\right. , \hspace{3mm} \left\{
\begin{array}{ccc}
J \Lambda_I J^{-1} &=&-\Lambda_I \\ J L_I J^{-1} &=&-L_I
\end{array}
\right. \eeq

we have the following generalized Kodaira identities: \beq
\label{ko1} \left\{
\begin{array}{ccc}
d^* &=&-[\Lambda_C, d_C] \\ d_C^* &=& [\Lambda_C,d]
\end{array}
\right. ,
\hspace{3mm}
\left\{
\begin{array}{ccc}
d &=&[L_C, d_C^*] \\ d_C &=&-[L_C,d^*]
\end{array}
\right. , \eeq \beq \left\{
\begin{array}{ccc}\label{ko3}
d_K^*&=&[ \Lambda_I,d_J ] \\ d_J^*&=&[d_K, \Lambda_I ]
\end{array}
\right. \eeq

Let $\varphi$ be a function with a compact support. With the help
of (\ref{ko1}),(\ref{ko3})  we easily derive the relation \bqa
(\Lambda_I)^2dd_Id_Jd_K(\varphi)=2d_I^*d_Id_K^*d_K(\varphi)=
2\Delta^2(\varphi) \eqa

For an arbitrary top degree differential form $\psi$ on a
 four dimensional K\"ahler manifold $M$  there is a simple relation:
$(\Lambda_C^2\psi) \hspace{0.6mm} vol_M=2 \psi$ where $vol_M$ is
the volume form on $M$. Therefore we have the formula: $$
dd_Id_Jd_K(\varphi)=vol_M \Delta^2 (\varphi). $$ $\blacksquare$

\vskip 4mm
\begin{center}
\large \bf Higher analytic hypertorsion forms
\end{center}
\vskip 4mm

The  conditions on the differential forms discussed in the
previous section naturally  arise when the characteristic classes
of vector bundles are considered. Suppose we have a vector bundle
$\mathcal{E}$ over the  Riemannian manifold $M$. According to
Chern-Weil theory the choice of the connection  on the bundle
allows to construct the Chern character with values in the closed
differential forms. The image in the cohomology lies in the
integer lattice $H^{even}(M,\mathbb{Z})$ and gives the topological
invariant of the bundle. The smooth deformations of the bundle do
not change the cohomology class of the corresponding differential
form and  the conditions (\ref{real1}),(\ref{real2}) are
satisfied. The exactness of the variation of the Chern form allows
to construct Chern-Simons differential forms. For instance
considering the second Chern class $c_2=-\frac{1}{8\pi^2} Tr
F\wedge F$ for a one dimensional family of the connections on the
bundle parameterized by the variable $t$  we get an example of the
Chern-Simons form $CS(A)=Tr (A\wedge dA +\frac{2}{3}A^3)$ through
the relation:

 \bqa
 \delta_t c_2(A)=-\frac{1}{8\pi^2}d \delta_t CS(A) \eqa

 In the case of  holomorphic bundles on the K\"ahler manifold
 the choice of a hermitian metric on the bundle  leads to
 the Chern character form subjected to the additional
 condition.
The corresponding cohomology classe  should be invariant under the
natural action of  $U(1)$ on the cohomology of the K\"ahler
manifolds. Thus we have all the conditions
(\ref{complex1}),(\ref{complex2}),(\ref{complex3}) satisfied and
this allows to define  Bott-Chern differential forms (i.e. see
\cite{BGS}for the detailed discussion) .

The next interesting case is a hyperholomorphic bundle on the
hyperk\"ahler manifold. Hyperholomorphic bundle is a hermitian
bundle vector bundle which is  holomorphic with respect to any of
the compatible holomorphic structures associated with the
hyperk\"{a}hler manifold. The corresponding characteristic classes
are subjected to the condition to be invariant with respect to the
action of the isotropy group on the cohomology \cite{Ver}. This
provides additional condition (\ref{quater3}) and allows to apply
the generalization of the $dd_c$-lemma from the first part of the
paper. Thus we have derived the existence of the fourth order
transgression of the Chern character form of an arbitrary
hyperholomorphic vector bundle. Note however that this arguments
is global and one could wonder if there exists a simple local
expression for the resulted differential form.

Below we give the explicit answer for one particular interesting
example. We consider the infinite dimensional hyperholomorphic
bundles naturally arising from the families of the hyperk\"ahler
manifolds supplied with a finite dimensional hyperholomorphic
bundle. We provide local construction of the fourth order
transgression in this case  and give the explicit formula for the
resulted resulted generalized higher torsion form (hypertorsion
form). Local construction for the general case of an arbitrary
hyperholomorphic bundle will be discussed elsewhere.

Consider the local universal family $\pi :M\times B \rightarrow B$
of the deformations of a hyperholomorphic bundle $\V$ with a
hermitian metric on the fiber $M$ parameterized by $B$. Let $\W$
be a corresponding universal bundle over $M\times B$. The family
of the Dirac operators $D=D^-+D^+$ acting along the fiber $M$ on
the twisted spinor bundles $\E=\V \otimes \s(M)=\V \otimes
\s(M)_+\oplus \V \otimes \s(M)_-$ defines the virtual index bundle
$\widetilde{\W} \equiv Ind(D^+)$ on the base of the fibration.
This provides two closed differential forms on the base $B$. The
first form is the product of the Chern character of $\W$ and
$\widehat{A}$ class of the tangent bundle to $M$ integrated over
the fiber of the projection. The other one  is the product of the
Chern class of $\widetilde{\W}$ supplied with the $L_2$ metric.
The local families index theorem of Atiyah and Singer \cite{AS}
claims that:
 \bqa \label{AS} ch(\widetilde{\W})=\pi_* [ch(\W)\widehat{A}(TM)]
\mbox{\, in \,}
 H^{even}(B,\Q) \eqa

One can construct a one-dimensional family of Quillen
superconnections acting in the associated infinite dimensional
hyperholomorphic bundle $\widehat{\W}$ of twisted spinor sections
$\Gamma(\E|_M,M)$ over $B$. This gives rise to the
representative of the Chern character in the differential forms
interpolating between the l.h.s. and r.h.s. of (\ref{AS}). Locally
over $B$ both parts of (\ref{AS}) are
 given by exact forms and by the general properties of the Chern classes of
hyperholomorphic bundles \cite{Ver} are $G_M=Sp(1)$-invarinat. We
derive the explicit formula for fourth order transgression of
their difference. This defines hypertorsion differential form for
the families of hyperholomorphic bundles. Let us start with short
description of  Quillen superconnection formalism. Consider
$\Z_2-$graded  vector bundle $\E$ and let $\tau$ be the operator
defining the $\Z_2$-grading on $\E$ i.e. $\tau=\pm 1$ on
$\E_{\pm}$. The algebra $End(\E)$ is naturally $\Z_2$ graded
algebra. We set a $\Z_2-$grading to the bundle of $\E$-valued
differential forms $\Lambda^*(\E,M)$ as a graded tensor product.
For $A , A' \in \Lambda^*(\E,M)$ the supercommutator $[A,A']$ is
given by: \bqa
 [A,A']=AA'-(-1)^{deg(A)\,deg(A')}A'A \eqa
$\Z_2$-grading allows to define supertrace $Str (A)$ as: \bqa
 Str(A)=Tr(\tau A)\eqa
with the property to be zero on supercommutators.

The form $Str(e^{-\nabla^2})$ is a closed differential form
representing  Chern character of the virtual bundle $\E_+\ominus
\E_-$ : \bqa
 ch(\E_+\ominus \E_-)=ch(\E_+)-ch(\E_-)=Str(e^{-\nabla^2}) \eqa

Thus defined Chern classes differ from the standard Chern classes
by the multiplication of the degree $2k$ components by $(2\pi
i)^k$. In the following we will always this normalization.

 Let us given an odd self adjoint operator V acting on $\E$
(i.e odd section of $End (\E)$). We could combine it with the
connection to get the Quillen superconnection on $\Z_2$-bundle
$\E$.
\begin{deff}
{A differential operator $\A: \Lambda^*(\E,M)\rightarrow
\Lambda^*(\E,M)$ of order $1$ with respect to  $\Z_2-$grading is a
Quillen's superconnection if $$ \A (\omega s)=(d\omega)
s+(-1)^{deg \, \omega}\omega \A(s) $$ }
\end{deff}
In fact, we could construct a family of the superconnections
depending on a real positive parameter $t\in \RR_+$: \bqa
 \A_t=\nabla+\sqrt{t}V \eqa
Here  $\nabla$ and $V$ are even and odd parts of the
superconnection.

The space $B$ of  local deformations of a hyperholomorphic bundle
$\V$ on $M$  is naturally supplied with a hyperk\"ahler structure.
We show that Quillen superconnections defined over the base $B$
and extended base $B\times \HH$ are hyperholomorphic.

\begin{propos}
{ Let $B$ be a space of local deformations of a  hyperholomorphic
vector bundle with a hermitian metric over a hyperk\"ahler
manifold $M$.

(i) The superconnection $\A_t=d^B +\sqrt{t}D$ is hyperholomorphic
over $B$.

(ii) Consider the operator $D_x=x^0 D+x^1 D_I +x^2 D_J +x^3 D_K$,
where $D_L=c(L)Dc(L)^{-1}$. Chose a hyperk\"ahler structure on
quaternionic plane  $x=x^0 +x^1 i+x^2 j+x^3 k\in \HH$  by
considering the right multiplication by quaternionic units $-i$,
$-j$,  $-k$.

Then the superconnection $\A_x=d^{\HH}+ d^B +D_x$ is
hyperholomorphic over $B\times\HH$.
 }
\end{propos}
\proof

(i) Let $C$ be a complex structure compatible with the
hyperk\"ahler structure on $B\times M$. We will use  the results
from the end of the Appendix B. The isomorphism $\E=\V\otimes \s
\cong \Lambda^{(*,0)}(\V,M)$,
$D=\sqrt{2}(\nabla_C'+(\nabla_C')^*)$ on $M$ leads to the
decomposition $D=D'+D''$ where $D'=\sqrt{2}\nabla'_C$
,$D''=\sqrt{2}(\nabla'_C)^*$ have the types $(1,0)$ and $(0,1)$.
Using the variant of the Kadaira identity:
$\{\overline{\partial}_C^B,(\nabla_C)^*\}=0$ we have: \bqa
\label{hypcon1}
 (\A_t'')^2=(\overline{\partial}_C^B+\sqrt{2t}(\nabla_C)^*)^2=0
 \eqa
 In particular
 for $\A_{t,C}=i(\A_t'-\A_t'')$ the following
 identities holds:
 \bqa \label{Aident}
 \A_{t,I}^2 = \A_t^2 \eqa

Note that here $C$ acts on the total tangent bundle to $B\times
M$.
 Taking into account that (\ref{hypcon1}) holds for any $C$ we
 infer that $\A_t$ is hyperholomorphic over $B$.

(ii) Note  that $[ad_{\phi}, D_x]=D_{\phi x}$, where
 $\phi=\phi^1I+\phi^2J+\phi^3 K$ is an arbitrary generator of $Sp(1)-$action
on twisted spinors. Let us start with the complex structure
defined by $I$. Then $\sqrt{t}D'$ is gauge equivalent to
$D_x'=z_1D'+\bar{z}_2D_J'$, and $\sqrt{t}D''$ is gauge equivalent
to $D_x''=\bar{z}_1D''+z_2D_J''$, where $$z_1=x^0+\sqrt{-1}x^1,
\hspace{3mm} z_2=x^2+\sqrt{-1}x^3 $$ $$t=|x|^2=(x_0)^2+ (x_1)^2+
(x_2)^2 +(x_3)^2 $$
 and $D_J'=JD'J^{-1}$
$D_J''=JD''J^{-1}$. The connection operators
$\partial^B+\sqrt{t}D'$ and $\bar{\partial}^B+\sqrt{t}D''$ are
gauge equivalent to $\partial^B +D_x'$ and $\bar{\partial}^B
+D_x''$. Therefore their squares are equal to zero.

Let us decompose $\A_x=\A'+\A''$, where $$ \A'=d\bar{z}_1
\frac{\partial}{\partial \bar{z}_1}+dz_2
\frac{\partial}{\partial z_2} +\partial^B+z_1D'+\bar{z}_2 D_J'
$$ $$ \A''=dz_1 \frac{\partial}{\partial z_1}+d\bar{z}_2
\frac{\partial}{\partial \bar{z}_2}+
\bar{\partial}^B+\bar{z}_1D''+z_2D_J'' $$
It is clear that
$(\A')^2=(\A'')^2=0$

Since $(\bar{z}_1, z_2)$ are holomorphic coordinates on $\HH$ with
respect to right multiplication by $-i$ the superconnection $\A$
is holomorphic on $B \times \HH$. The same arguments work for any
compatible complex structure $C$. Taking into account the
isomorphism $\S\otimes\V \simeq\Lambda_C^{(*,0)}(\V,M)$ for any
complex structure $C$ we conclude that $\A_x$ is hyperholomorphic
over $B \times \HH$. $\blacksquare$


Consider the Chern character form defined by superconnection
$\A_t$ (see \cite{B}, \cite{BGV}, \cite{BGS}):

\bqa ch(\A_t)=Stre^{-\A_t^2} \eqa  It interpolates between Chern
character form $ch(\widetilde{\W})$ for the $L_2$-metric on the
index bundle $\widetilde{\W}$ and characteristic class
$ch(\W)\widehat{A}(TM)$ integrated along the fiber:

\bqa ch(\widetilde{\W})=\lim\limits_{t\rightarrow\infty}Stre^{-\A_t^2}\\
\int\limits_{M}ch(\W)\widehat{A}(TM)=\lim\limits_{t\rightarrow
0}Stre^{-\A_t^2} \eqa
 where $ch(\W)$ is
Chern character form of the canonical connection on the bundle
$\W$ over $M\times B$ and $\widehat{A}(TM)$ is  the multiplicative
genus given by the power series:
 $$ \widehat{A}(x)=\frac{\frac{x}{2}}{sinh(\frac{x}{2})}$$
of the curvature of the Levi-Chivita connection over $M$.

\begin{theorem}\label{th1}
{ The following transgression formula holds: \bqa \label{tfirst}
ch(\widetilde{\W})-\int\limits_{M}ch(\W)\widehat{A}(TM)=\frac{1}{24} d^B
d^B_I d^B_J d^B_K \beta \eqa where \bqa \label{tsecond}
\beta=\sum_{C=I,J,K}\int\limits_{\rightarrow 0}^{+\infty} Str
(\int\limits_0^1 ad_C e^{-\tau\A_t^2} ad_C e^{-(1-\tau)\A_t^2}
d\tau ) \frac{dt}{t} \eqa
 is  a higher hypertorsion differential form. The zero degree part of  $\beta$
 $$
\beta_0=\int\limits_{\rightarrow 0}^{+\infty} Str(
\sum_{C=I,J,K}(ad_C)^2 e^{-tD^2}\frac{dt}{t})$$ may be expressed
in terms of the Laplace operators $\Delta_q$ acting on $q$-froms:
$D^2 =\sum_q \Delta_q $ as the logarithm of the hypertorsion
$T_h$:
 \bqa \beta_0 &=& 3log T_h  \\
  T_h &=& \prod\limits_{q=0}^{2k}(det'\Delta_q)^{(-1)^q q^2} \eqa
 }
\end{theorem}

The definition of the "regularized" integral
$\int\limits_{\rightarrow 0}^{\infty} f(t)\frac{dt}{t}$ is given
in Appendix A.

First we give a simple local argument in favor of the existence of
the fourth order transgression  and then give the formal proof of
the theorem.

Taking into account the identity which follows from Proposition
\ref{lapl} :

  \bqa d^{\HH} d_I^{\HH} d_J^{\HH}
d_K^{\HH} \ln |x|^2=16\pi^2 (\delta_{\infty}-\delta_0) \eqa we
have the following representation:  \bqa Str
e^{-\A_t^2}|^{+\infty}_0 =\frac{1}{16\pi^2}\int\limits_{\HH} Str
e^{-\A_x^2} d^{\HH} d_I^{\HH} d_J^{\HH} d_K^{\HH} ln|x|^2 \eqa

 Since
$\A_x=d^{\HH}+d^B +D_x$ is hyperholomorphic over $B\times \HH$, we
obtain  \bqa Str e^{-\A_t^2}|^{+\infty}_0=\frac{1}{16\pi^2}d^B
d^B_I d^B_J d^B_K
 \int\limits_{\HH} Str e^{-\A_x^2} ln|x|^2
\eqa

This leads to the fourth-order transgression of the difference of
the Chern character forms: \bqa
 ch(\widetilde{\W})
-\int\limits_{M}ch(\W)\widehat{A}(M)=\frac{1}{16\pi^2}d^B d^B_I
d^B_J d^B_K
 \int\limits_{\HH} Str e^{-\A_x^2} ln|x|^2
\eqa

This representation provides the  direct generalization of the
representation for the higher holomorphic torsion form in terms of
the integration over auxiliary complex plane given  \cite{GS}.

One could reduce the expression in r.h.s. to the one given in the
Theorem \ref{th1}. However to make analytic regularization more
explicit we proceed with the direct derivation of
(\ref{tfirst})(\ref{tsecond}).

\vspace{5mm} \centerline{\bf Proof of the theorem.} \vspace{5mm}
Let us start with the following lemma.


\begin{lemma}
We have: \beq t\pr_t (t\pr_t +1)Str e^{-\A_t^2}= \frac{1}{8}d^B
d_I^B d_J^B d_K^B Str (\int\limits_0^1  ad_I e^{-\tau\A_t^2} ad_I
e^{-(1-\tau)\A_t^2} d\tau ) \label{useful} \eeq
\end{lemma}
\proof

By simple calculation using the relations $\A_t^2=\A_{t,I}^2$ ,
$[\A_{t,I},\A_t^2]=0$ (see \cite{BGS} for  similar considerations)
we get: $$ \pr_t Str e^{-\A_t^2} = -\frac{1}{2\sqrt{t}} Str
([\A_t, D] e^{-\A_t^2}) = -d^B \frac{1}{2\sqrt{t}} Str (D
e^{-\A_t^2})= $$ $$ =-d^B \frac{1}{2 t} Str ([\A_{t,I}, ad_I]
e^{-\A_t^2})= -\frac{1}{2t} d^B d_I^B Str (ad_I e^{-\A_t^2}) $$

Thus we have the relation: \beq t\pr_t Str e^{-\A_t^2}=
-\frac{1}{2} d^B d_I^B Str (ad_I e^{-\A_t^2}) \label{first} \eeq

Applying the formula of differentiation:  \bqa
 \frac{d}{dt}e^{-A(t)}=-\int_0^1 e^{-
 \tau A(t)}\frac{d}{dt}A(t)e^{-(1-\tau)
 A(t)}d\tau \eqa
we derive  $$ \pr_t Str (ad_I e^{-\A_t^2})= -\frac{1}{2\sqrt{t}}
Str (ad_I \int\limits_0^1 e^{-\tau\A_{t,J}^2} [\A_{t,J},D_J]
e^{-(1-\tau)\A_{t,J}^2} d\tau )= $$ $$ =\frac{1}{2\sqrt{t}}(
-d_J^B Str (ad_I \int\limits_0^1 e^{-\tau\A_{t}^2} D_J
e^{-(1-\tau)\A_{t}^2} d\tau )+
 Str ([\A_{t,J},ad_I]\int\limits_0^1
e^{-\tau\A_{t}^2} D_J e^{-(1-\tau)\A_{t}^2} d\tau ))= $$ $$
=-\frac{1}{2t}d_J^B Str (ad_I \int\limits_0^1 e^{-\tau\A_{t}^2}
[\A_{t,K},ad_I] e^{-(1-\tau)\A_{t}^2} d\tau )- \frac{1}{2} Str
(D_K \int\limits_0^1 e^{-\tau\A_{t}^2} D_J e^{-(1-\tau)\A_{t}^2}
d\tau ) $$
The first part of the expression is equal to
$-\frac{1}{4t}d_J^B d_K^B Str (\int\limits_0^1 ad_I
e^{-\tau\A_{t}^2} ad I e^{-(1-\tau)\A_{t}^2} d\tau )$. Acting by
$t\pr_t$
on (\ref{first}) and using the result of the previous
calculation we get \beq (t\pr_t )^2 Str e^{-\A_t^2}=
\frac{1}{8}d^B d_I^B d_J^B d_K^B Str (\int\limits_0^1 ad_I
e^{-\tau\A_t^2} ad_I e^{-(1-\tau)\A_t^2} d\tau )+ \frac{t}{4} d^B
\alpha, \label{second} \eeq
 where
$$ \alpha=d^B_I Str (D_K \int\limits_0^1 e^{-\tau\A_{t}^2} D_J
e^{-(1-\tau)\A_{t}^2} d\tau ) $$

The next step is to obtain $\alpha$. $$ \alpha=\int\limits_0^1
d\tau (Str ( [\A_{t,I},D_K] e^{-\tau\A_{t}^2} D_J
e^{-(1-\tau)\A_{t}^2}- D_K e^{-\tau\A_{t}^2} [\A_{t,I}, D_J]
e^{-(1-\tau)\A_{t}^2} )) $$ Note that using: \bqa \label{commut}
 [\A_{t,I},ad_I]=\sqrt{t}D \eqa
one could get the following relations:\bqa
[\A_{t,I},D_K]=\frac{1}{\sqrt{t}}[\A_{t,I},[\A_{t,I},ad_J]]=
\frac{1}{\sqrt{t}}[\A_t^2, ad_J]\eqa \bqa
[\A_{t,I},D_J]=-\frac{1}{\sqrt{t}}[\A_{t,I},[\A_{t,I},ad_K]]=-
\frac{1}{\sqrt{t}}[\A_t^2, ad_K] \eqa
 So we have $$
\alpha=\frac{1}{\sqrt{t}}\int\limits_0^1 d\tau Str ([\A_t^2, ad_J]
e^{-\tau\A_{t}^2} D_J  + D_K e^{-\tau\A_{t}^2}[\A_t^2, ad_K])
e^{-(1-\tau)\A_{t}^2}= $$ $$ =\frac{1}{\sqrt{t}} \int\limits_0^1
d\tau
 Str (ad_J \pr_{\tau}
e^{-\tau\A_{t}^2} D_J e^{-(1-\tau)\A_{t}^2}- D_K \pr_{\tau}
e^{-\tau\A_{t}^2} ad_K e^{-(1-\tau)\A_{t}^2})= $$ $$
=\frac{1}{\sqrt{t}} Str (ad_J [e^{-\A_t^2}, D_J] -D_K
[e^{-\A_t^2}, ad_K])=\frac{1}{\sqrt{t}} Str e^{-\A_t^2} ([D_J,
ad_J] + $$ $$ +[D_K, ad_K])=\frac{2}{\sqrt{t}}Str D e^{-\A_t^2}.
$$

Taking into account (\ref{commut}) we derive: \beq
\alpha=\frac{2}{t} d^B_I Str (ad_I e^{-\A_t^2}). \label{third}
\eeq

After substitution of (\ref{third}) in (\ref{second}) we obtain:
\beq (t\pr_t )^2 Str e^{-\A_t^2}= \frac{1}{8}d^B d_I^B d_J^B d_K^B
\Phi (t)+
 \frac{1}{2} d^B d^B_I Str (ad_I e^{-\A_t^2}),
\eeq where $$ \Phi (\kappa )=Str (\int\limits_0^1  ad_I
e^{-\tau\A_{\kappa}^2} ad_I e^{-(1-\tau)\A_{\kappa}^2} d\tau) $$

By using (\ref{first}) we immediately prove the lemma.
$\blacksquare$

Let us apply (\ref{App1}) to  the combination
$G_{\A}(t)=Stre^{-\A^2_t}-ch(\widetilde{\W})$. Note that
$G_{\A}(t)=O(t^{-\frac{1}{2}})$ as $t\rightarrow\infty$. Thus we
have: \bqa
 \int^{\infty}_{\rightarrow 0}(t\pr_t(t\pr_t+1)(Stre^{-\A^2_t}
 -ch(\widetilde{\W})\frac{dt}{t}= \nonumber \\
 = -(t\pr_t+1)G_{\A}|_{t=0}=ch(\widetilde{\W})-Stre^{-\A^2_t}(0)\eqa
Taking into account (\ref{useful}) we have proved the first part
of the theorem.

Now let us prove the formula for the zero-degree part of
hypertorsion form.

It is clear that $$ \beta_0=3\int\limits_{\rightarrow 0}^{+\infty}
Str (\int\limits_{0}^1  ad_I e^{-\tau tD^2} ad_I e^{-(1-\tau)t
D^2} d\tau ) \frac{dt}{t}=3\int\limits_{\rightarrow 0}^{+\infty}
Str (ad_I^2 e^{-tD^2})\frac{dt}{t}$$

The following identity from the Appendix A being applied to the
trace of the positive self-adjoint operator $\hat{H}$: $$
\int\limits_{\rightarrow 0}^{\infty}
Tre^{-t\hat{H}}\frac{dt}{t}=-log det'\hat{H}, $$ immediately leads
to the representation of the zero-degree part of the hypertorsion
form $\beta_0$ in terms of infinite determinants.

Recall that $ad_I$ acts on $(q,0)-$forms on $4k$-dimensional
K\"ahler manifold as $i(q-k)$. So $(ad_I)^2=-(q-k)^2$ and hence $$
\beta_0 \equiv 3\log T_h=3\sum\limits_{q=0}^{2k}(-1)^q (q-k)^2 Tr
\Delta_q^{-s} =3\sum\limits_{q=0}^{2k}(-1)^q (q-k)^2 log
det'\Delta_q, $$ where $\Delta_q$ is the Laplace operator acting
on $\W$ valued $(q,0)-$forms, $\Delta=D^2=\sum_{q}\Delta_q$, and
$det'\Delta_q$ is the regularalized determinant. We have $$
T_h=\prod\limits_{q=0}^{2k}det'\Delta_q^{(-1)^q (q-k)^2}. $$

 The usual analytic torsion for holomorphic bundle over a complex
 manifold (see \cite{BGV}, \cite{BGS}, \cite{RS})
 is given by $T=\prod\limits_{q=0}^{2k}(det'\Delta_q)^{q(-1)^q}$.
In the following lemma we prove that analytic torsion for a
hyperholomorphic bundle over a hyperk\"ahler manifold is trivial.


\begin{lemma}
{Let $M$ be a hyperk\"ahler manifold and let $\V$ be a
hyperholomorphic bundle over $M$. Then $T=1$.}
\end{lemma}
\proof

Let us write down the expression for $log T$ in the following
form: $$ \log T= \sum\limits_{q=0}^{2k}(-1)^q
(-\frac{\partial}{\partial s}_{|s=0})tr(h+k)\Delta_q^{-s}, $$
where $h=\frac{1}{i}ad_I=q-k$. Let $S^{\pm}_{\lambda}$ be the
eigen-spaces of $\Delta_{\pm}$, where $$\Delta_+
=\sum_{q-even}\Delta_q, \hspace{4mm} \Delta_-
=\sum_{q-odd}\Delta_q.$$ Since $D\Delta_{\pm}=\Delta_{\mp}D$ then
$D:S^{\pm}_{\lambda} \simeq S^{\mp}_{\lambda}$. Thus  we conclude
that $$\sum\limits_{q=0}^{2k}(-1)^q tr\Delta_q^{-s}=0.$$

Moreover, the Laplace operator is $Sp(1)-$invariant, so the
eigen-subspaces $S^{\pm}_{\lambda}$ are $\sp (1)-$modules. It
follows that $tr(h)_{|S^{\pm}_{\lambda}}=0$ and
$trh\Delta^s_+=trh\Delta^s_-=0$. Therefore $\ln T=0$ and $T=1$.
$\blacksquare$

This lemma implies that $$
T_h=\prod\limits_{q=0}^{2k}det'\Delta_q^{(-1)^q q^2}. $$
$\blacksquare$

There is an interesting particular case of the theorem we have
proved. Let $dim M=4$ and $ch_{[2]}$ be a component of the Chern
character taking values in four-forms . In this case
$T_h=(det'\Delta_0)^2$ and we have: $$
ch_{[2]}(\widetilde{\W})=(\int\limits_{M}ch(\W)\widehat{A}(M))_{[2]}+\frac{1}{4}d^B
d_I^B d_J^B d_K^B log (det'\Delta_0) $$ If in addition $dimB=4$
then $$
ch_{[2]}(\widetilde{\W})=(\int\limits_{M}ch(\W)\widehat{A}(M))_{[2]}+
\frac{1}{4} vol_B \Delta_B^2 log (det'\Delta_0) $$ This formula
was proposed in the physical literature in \cite{Sch},\cite{BB}
for the case of $M=T^4$ and in \cite{CG} (see also \cite{BFST})
for the instantons over $M=\R^4$.

\vskip 4mm

\begin{center}
\large \bf Appendix A: Regularization of integrals.
\end{center}

In this appendix we define the regularization of some class  of
the integrals using analytic continuation. This regualrization is
a standard tool in the theory of higher analytic torsion
\cite{B,BGS}

Let $G(t)$ be a continuous function defined for $t>0$ with
sufficiently rapid decay as $t \rightarrow \infty$. We also assume
that it has an asymptotic expansion as $t\rightarrow 0$:
 \bqa \label{assymp} G(t)=\sum_{i=-n}^{0}G_it^i+O(t) \eqa

The following integral: \bqa
\zeta_G(s)=\frac{1}{\Gamma(s)}\int^{\infty}_0 G(t)t^{s-1}dt \eqa
converges for $Re(s)>n$ and  has the analytic extension to the
whole complex plane. We define the value of ( in general
divergent) integral as follows: \bqa
 \int^{\infty}_{\rightarrow 0} G(t)\frac{dt}{t}\equiv \zeta_G'(0) \eqa
For instance by this definition for $G(t)=e^{-th}$ we have: \bqa
 \int^{\infty}_{\rightarrow 0}e^{-th} \frac{dt}{t}=\zeta_{e^{-th}}'(0)=-log(h)
 \eqa
Note that thus defined integral has the usual property for the
total derivative  of a regular function:
 \bqa \label{fulder}
 \int^{\infty}_{\rightarrow 0}(t\partial_t F(t)) \frac{dt}{t}=-F(0)
 \eqa
If $F(t)$ has a more general behaviour (\ref{assymp}), the value
at $t=0$ of the regular part of $F(t)$ appears in r.h.s of
(\ref{fulder}) instead of $F(0)$   .

 We need the following consequence of this property. Consider
the
 regularized integral of the function
$H(t)=t\partial_t(t\partial_t+1) G(t)$ for regular $G$:  \bqa
 \zeta_H(s)= \frac{1}{\Gamma(s)}\int^{\infty}_0 (t\partial_t(t\partial_t+1) G(t))t^{s-1}dt \eqa

Then the following obvious identity holds: \bqa \zeta_H'(0)=-G(0)
\eqa and we have: \bqa \label{App1}
 \int^{\infty}_{\rightarrow 0}
(t\partial_t(t\partial_t+1) G(t))\frac{dt}{t}=-G(0) \eqa

\vskip 4mm

\begin{center}
\large \bf Appendix B: Dirac operator on  hyperk\"ahler manifolds
\end{center}

Here we recall the interrelation of spin structures and complex
structures  on the K\"ahler and hyperk\"ahler manifolds with the
emphasis on the properties of the Dirac operator.

\vskip 3mm
\begin{deff}\cite{BGV}
{Let $V$ be a real $n-$dimensional vector space with positive
quadratic form  $g$. The Clifford algebra of $(V,g)$, denoted by
$Cl(V)$, is the algebra over $\R$ generated by $V$ with the
relations $xy+yx=-2g(x,y)$. A self-adjoint Hermitian module $E$ of
$Cl(V)$ is called a Clifford module.}
\end{deff}

Let $e^i$, $i=1,...n$ be an orthogonal basis of $V$ and let $c^i$
be an element of $Cl(V)$ corresponding to $e^i$. One can extend
this map to the isomorphism of graded $O(V)$ modules $c: \Lambda^*
(V) \rightarrow Cl(V)$ by sending $e^{i_1}\wedge ... \wedge
e^{i_k}\mapsto c^{i_1}....c^{i_k}$. The chirality operator $\Gamma
=(\sqrt{-1})^{[\frac{n+1}{2}]} c^1...c^n$ satisfies $\Gamma^2=1$
and defines a $\Z_2-$grading on $Cl(V)\otimes_{\R}\C$. Taking
$v\in V$ one can check that $c(v)\Gamma=(-1)^{n+1}\Gamma c(v)$.

The subspace $Cl^2 (V)=c(\Lambda^2(V))$  is a Lie subalgebra of
$Cl(V)$ which is isomorphic to $\so (V)$ under the map $\tau :
Cl^2 (V) \simeq \so (V)$, $\tau (a)x=[a,x]$, $v\in V$. The group
$Spin(V)$ is  obtained by exponentiation of the Lie algebra
$Cl^2(V)$ inside the Clifford algebra $Cl(V)$.

Let $V$ be a hermitian vector space with a complex structure $I$.
Let $\omega_I=g(I\cdot,\cdot)$ be the corresponding real
nondegenerate $2-$form, which can be considered by duality as an
element of $\Lambda^2(V)$. Using the isomorphism $Cl^2 (V) \simeq
\so (V)$ it is possible to show that $c(\omega^I)=2I$.

Let us decompose $V^c=V\otimes_{\R}\C$ into holomorphic and
antiholomorphic parts $V^c=W\oplus \bar{W}$. Since $W$ and
$(\bar{W})$ are isotropic subspaces with respect to the scalar
product extended by complex linearity, then $Cl(W)$ and
$Cl(\bar{W})$ are commutative graded algebras.

Let us define an irreducible Clifford module, denoted by $S$,
which is called a spin module as a hermitian complex space
$\Lambda^*(W)$, provided with the following Clifford action $$
c(w)= \left\{ \be{cc} \sqrt{2}\epsilon(w), & w\in W, \\
-\sqrt{2}\iota (\bar{w}), & w\in \bar{W}, \ee \right. $$ where
$\epsilon(w)$ is the exterior product of $w$, and $\iota(\bar{w})$
is the contraction with the hermitian dual covector. The spin
representation constructed as above has a unique up to
multiplication by unitary complex numbers normalized vacuum vector
$|1\rangle$,
 which satisfies  the following conditions:
$ Cl(\bar{W})|1\rangle = 0$ and $ S = Cl(W)|1\rangle$. The module
$S$ descends a natural $\Z-$grading from the space
$Cl(W)=\bigoplus_{q=0}^n Cl^q(W)$. One can verify that $$
c(\omega^I)_{|Cl^q(W)|1\rangle}=i(2q-n). $$

Let $V$ be a self-adjoint $\HH -$module of dimension $4k$.
  The correspondence between the families of complex
structures $C$ and the associated $2-$forms $\omega^C$ immediately
leads to the inclusion $\sp(1)\hookrightarrow Cl(V)$,
$C\rightarrow \frac{1}{2}c(\omega^C)$, which can be exponentiated
inside $Cl(V)$ to the inclusion  $Sp(1)\hookrightarrow Spin(V)$.
Therefore the group of unitary quaternions $Sp(1)$ acts on the
spin module $S$, such that $c(x)c(v)c(x)^{-1}s=c(x(v))s$, where
$x\in Sp(1)$, $v\in V^c$, $s\in S$. If $dimV=4$ then the subspace
of self-dual $2-$forms $\Lambda^2_+ (V)$ is spanned by $\omega^C$
(\cite{AHS}).

Since $\sp(1)\otimes_{\R}\C=\sll(2,\C)$, one can choose the
$\sll_2-$generators $h=\frac{1}{2i}c(\omega^I)$,
$e=\frac{1}{4}(c(\omega^J)-ic(\omega^K))$,
$f=-\frac{1}{4}(c(\omega^J)+ic(\omega^K))$ with the relations
$[h,e]=-2f$, $[h,f]=-2f$, $[e,f]=h$.

Let us consider the spin module $S$ over $Cl(V)$ as the linear
space of $(p,0)-$forms with respect to $I$. Then we have the
following simple property:
\begin{propos}\label{forms}
{The form $\Omega =\frac{1}{4}(\omega^J-i\omega^K)$ is of type
$(2,0)$ with respect to $I$. Moreover, the operator $e$ acts as
the exterior product with $\Omega$ and the operator $f$ acts as
the contraction with $\Omega$. }
\end{propos}

Since
$h_{|\Lambda^{q,0}(V)}=q-k$, where $q=0,...2k$ we see that the
operator $h$ defines $\Z-$grading on $S$.

\begin{deff}\cite{BGV}
{A Clifford module $\E$ over an even dimensional Riemannian
manifold $M$ is a $\Z_2-$graded hermitian bundle of Clifford
modules $\E=\E_+ \oplus\E_-$ over a bundle $Cl(M)$ of Clifford
algebras with the unitary connection $\nabla^{\E}$, which is
compatible with the Levi-Chivita connection, extended to $Cl(V)$.
If $\V$ is a hermitian vector bundle, then $\V\otimes\E$ is the
twisted Clifford module with the Clifford action $1\otimes c(a)$
and with the connection $\nabla^{\V \otimes\E}=\nabla^{\V}\otimes
1+1\otimes\nabla^{\E}.$}
\end{deff}

The action of one-forms $\Lambda^1(M)$ on $\E$ defines a
$C^{\infty}(M)-$linear
 morphism of bundles over $M$, written as $\Lambda^1(\E_{\mp}, M)\rightarrow \Gamma
(\E_{\pm},M)$.
 So one can introduce a generalized Dirac operator,
 acting  as follows:
$$ D:
\Gamma(\E_{\mp},M)\stackrel{\nabla^{\E_{\mp}}}{\longrightarrow}
\Lambda^1(\E_{\mp},M)\rightarrow \Gamma (\E_{\pm},M). $$

Let us define a {\sl spin bundle} as a Clifford bundle of spin
models. Given an almost complex structure $I$ one can construct a
bundle of spin modules. If the Riemannian manifold is K\"ahler we
have a subbundle  of $(p,0)-$forms $\Lambda^{*,0}(M)$. This bundle
has the structure of the Clifford module. More generally there is
\begin{propos}\cite{BGV}
{Let $\V$ be a holomorphic vector bundle with hermitian metric on
a K\"ahler manifold. The tensor product of the Levi-Chivita
connection on $\Lambda^{*,0}(M)$ with the canonical connection on
 $\V$ gives  Clifford connection on the Clifford module
$\Lambda^{*,0}({\cal V},M)$. Let $\nabla'$ be  (1,0)-part of the
connection  on the bundle $\Lambda^{*,0}({\cal V},M)$. Then the
Dirac operator on the corresponding Clifford module is
$\sqrt{2}(\nabla' +\nabla'^*)$. }
\end{propos}

Let $M$ be a hyperk\"{a}hler manifold. Then there is a covariantly
constant inclusion of $Sp(1)$ as the gauge subgroup of Clifford
bundle's sections. Using a fixed complex structure $I$ from the
hyperk\"{a}hler family of complex structures, we can construct a
spin bundle $\s$ over $M$ as above. On the hyperk\"ahler manifold
thus constructed spin bundle does not actually depend on the
choice of the complex structure $I$. This observation may be
exploited to prove the following proposition.
\begin{propos}
{Let $\V$ be a hyperholomorphic vector bundle with hermitian
metric on a hyperk\"{a}hler manifold. Then for any compatible
complex structure $C$  twisted spinor bundle
$\V\otimes\s=\V\otimes\s^+ \oplus \V\otimes\s^-$ is isomorphic to
$\Lambda_C^{*,0}({\cal V},M)=\Lambda_C^{even,0}({\cal V},M) \oplus
\Lambda_C^{odd,0}({\cal V},M)$. Under this isomorphism the Dirac
operator goes into $\sqrt{2}(\nabla'_C + (\nabla'_C)^*)$. }
\end{propos}

Consider the spin bundle $\s$ over $Cl(M)$ as the bundle of
$(*,0)-$forms with respect to $I$. Then $\Omega
=\frac{1}{4}(\omega^J-i\omega^K)$ is covariantly constant
$(2,0)-$form, therefore $\Omega$ is holomorphic. As a direct
consequence of this fact and the Proposition \ref{forms} one can
obtain, that the operator $e$ acts as the exterior product with
$\Omega$ and the operator $f$ acts as the contraction with
$\Omega$ on the space of twisted spinors. The last space is
identified with the space of $(*,0)-$forms.


\end{document}